\documentclass[12pt,reqno]{amsart}
\usepackage{amsmath,amssymb,amscd,latexsym,amsthm,mathrsfs}
\usepackage[unicode]{hyperref}
\usepackage{hypbmsec}
\ifx\pdfoutput\undefined\else
\hypersetup{
colorlinks=true,
linkcolor=blue,
citecolor=blue,
urlcolor=blue,
filecolor=blue,
bookmarksnumbered=true,
pdfstartview=FitH,
pdfhighlight=/N
}
\fi
\newtheorem*{maintheorem}{Main theorem}

\newtheorem{lemma}{Lemma}

\theoremstyle{remark}

\renewcommand\pmod[1]{\;(\operatorname{mod}#1)}

\begin{document}

\hypersetup{pdfauthor={Jes\'us Guillera, Wadim Zudilin},%
pdftitle={``Divergent'' Ramanujan-type supercongruences}}

\title{``Divergent'' Ramanujan-type supercongruences}

\author{Jes\'us Guillera}
\address{Av.\ Ces\'areo Alierta, 31 esc.~izda 4$^\circ$--A, Zaragoza, SPAIN}
\email{jguillera@gmail.com}

\author{Wadim Zudilin}
\address{School of Mathematical and Physical Sciences,
The University of Newcastle, Callaghan NSW 2308, AUSTRALIA}
\email{wadim.zudilin@newcastle.edu.au}

\date{12 July 2010}

\begin{abstract}
``Divergent'' Ramanujan-type series for $1/\pi$ and $1/\pi^2$
provide us with new nice examples of supercongruences
of the same kind as those related to the convergent cases.
In this paper we manage to prove three of the supercongruences
by means of the Wilf--Zeilberger algorithmic technique.
\end{abstract}

\maketitle

\section{Introduction}
\label{s1}

The following two supercongruences in the spirit of~\cite{Zu},
\begin{align}
\sum_{n=0}^{p-1}\frac{(\frac12)_n^3}{n!^3}(3n+1)2^{2n}
&\equiv p\pmod{p^3} \quad\text{for $p>2$},
\label{J1a}
\\
\sum_{n=0}^{p-1}\frac{(\frac12)_n^5}{n!^5}(10n^2+6n+1)(-1)^n2^{2n}
&\equiv p^2\pmod{p^5} \quad\text{for $p>3$},
\label{J2a}
\end{align}
correspond to divergent Ramanujan-type series for $1/\pi$ and $1/\pi^2$, respectively
(cf.\ Section~\ref{sdiv} below); the letter~$p$ is reserved for primes throughout the paper.
Furthermore, we have more congruences of this kind:
\begin{align}
\sum_{n=0}^{p-1}\frac{(\frac12)_n^3}{(1)_n^3}(3n+1)(-1)^n2^{3n}
&\equiv \biggl(\frac{-1}{p}\biggr)p\pmod{p^3} \quad\text{for $p>2$},
\label{3F2-zu3}
\displaybreak[2]\\
\sum_{n=0}^{p-1}\frac{(\frac12)_n^3}{(1)_n^3}(21n+8)2^{6n}
&\equiv 8p\pmod{p^3} \quad\text{for $p>2$},
\label{3F2-zu2}
\displaybreak[2]\\
\sum_{n=0}^{p-1}\frac{(\frac12)_n(\frac14)_n(\frac34)_n}{(1)_n^3}(5n+1)(-1)^n\biggl(\frac43\biggr)^{2n}
&\overset?\equiv \biggl(\frac{-3}{p}\biggr)p\pmod{p^3} \quad\text{for $p>3$},
\label{3F2-zu4}
\displaybreak[2]\\
\sum_{n=0}^{p-1} \frac{(\frac12)_n(\frac14)_n(\frac34)_n}{(1)_n^3}
(35n+8)\biggl(\frac43\biggr)^{4n}
&\overset?\equiv 8p \pmod{p^3} \quad\text{for $p>3$},
\label{3F2-zun}
\displaybreak[2]\\
\sum_{n=0}^{p-1}\frac{(\frac12)_n(\frac13)_n(\frac23)_n}{(1)_n^3}(11n+3)\biggl(\frac{27}{16}\biggr)^n
&\overset?\equiv 3p\pmod{p^3} \quad\text{for $p>2$},
\label{3F2-zu5}
\displaybreak[2]\\
\sum_{n=0}^{p-1}\frac{(\frac12)_n^5}{(1)_n^5}(205n^2+160n+32)(-1)^n2^{10n}
&\overset?\equiv 32p^2\pmod{p^5} \quad\text{for $p>3$},
\label{5F4-zu2}
\displaybreak[2]\\
\nonumber
\sum_{n=0}^{p-1}\frac{(\frac12)_n(\frac13)_n(\frac23)_n(\frac14)_n(\frac34)_n}{(1)_n^5}
\qquad\qquad\qquad \\[-7.5pt] \times
(172n^2+75n+9)(-1)^n\biggl(\frac{27}{16}\biggr)^n
&\overset?\equiv 9p^2\pmod{p^5} \quad\text{for $p>2$}.
\label{5F4-zu4}
\end{align}
Here $\bigl(\frac{\cdot}p\bigr)$ is the Legendre symbol and the Pochhammer notation
$(a)_b$ is used for denoting $\Gamma(a+b)/\Gamma(b)$ also in the cases when $b$~is not
a non-negative integer; of course, if $b=n\in\mathbb Z_{\ge0}$ we have, as usual,
$(a)_n=\prod_{k=0}^{n-1}(a+k)$ with the convention that the empty product equals~1.
The question mark indicates that the corresponding supercongruence remains conjectural;
the non-questioned entries~\eqref{J1a}--\eqref{3F2-zu3} are proved in this paper
by extending the method of~\cite{Zu}, while the supercongruence~\eqref{3F2-zu2}
(even in a more general form) is shown by Zhi-Wei Sun in his preprint~\cite{Su}.

Note that we can sum in \eqref{J1a}, \eqref{J2a}, \eqref{3F2-zu2}, \eqref{3F2-zu3},
and \eqref{5F4-zu2} up to $\frac{p-1}2$,
since the $p$-adic order of $(\frac12)_n/n!$ is~1 for $n=\frac{p+1}2,\dots,p-1$.

\begin{maintheorem}
The following supercongruences take place\textup:
\begin{align}
\sum_{n=0}^{(p-1)/2}\frac{(\frac12)_n^3}{n!^3}(3n+1)2^{2n}
&\equiv p\pmod{p^3} \quad\text{for $p>2$},
\label{J1}
\\
\sum_{n=0}^{(p-1)/2}\frac{(\frac12)_n^5}{n!^5}(10n^2+6n+1)(-1)^n2^{2n}
&\equiv p^2\pmod{p^5} \quad\text{for $p>3$},
\label{J2}
\\
\sum_{n=0}^{p-1}\frac{(\frac12)_n^3}{(1)_n^3}(3n+1)(-1)^n2^{3n}
&\equiv(-1)^{(p-1)/2}p\pmod{p^3} \quad\text{for $p>2$}.
\label{J4}
\end{align}
\end{maintheorem}

We find quite illogical that our strategy based on the creative Wilf--Zeilberger
theory~\cite{WZ} of WZ-pairs allows us to do only three entries from the list \eqref{J1a}--\eqref{5F4-zu4};
a very similar lack of luck was reported in~\cite{Zu}. Although we have
WZ-pairs for \eqref{3F2-zu2}--\eqref{5F4-zu2} as well, they seem to be
quite helpless for showing the corresponding congruences modulo the expected
powers of~$p$. Because of the clear relationship of such congruences
with Ramanujan's formulae for $1/\pi$ and their generalizations
(see \cite{Zu} and Section~\ref{sdiv}),
we do expect a more universal method for proving the Ramanujan-type supercongruences.

In Section~\ref{s2} we present auxiliary congruences, some of them are remarkable
in their own. Section~\ref{s3} contains the proofs of~\eqref{J1}--\eqref{J4}.
The final Section~\ref{sdiv} reviews the ``divergent'' Ramanujan-type series
for $1/\pi$ and $1/\pi^2$ as our motivation to the the above list of supercongruences.

\medskip
After posting the preprint online we were informed by Zhi-Wei Sun that he had
experimentally and independently discovered the congruences~\eqref{J1a}--\eqref{5F4-zu2},
but also proved in~\cite{Su} a more general than in~\eqref{3F2-zu2} supercongruence.
We thank him for attracting our attention to his work.

\section{Precongruences}
\label{s2}

In this section we summarize our needs for proving the supercongruences
of the main theorem.

\begin{lemma}
\label{standard}
The following congruences are valid\textup:
\begin{align}
\sum_{n=1}^{(p-1)/2}\frac{\binom{2n}n}n
&\equiv0\pmod p \quad\text{for $p>3$},
\label{st1}
\\
\sum_{n=1}^{(p-1)/2}\frac{(-1)^n\binom{2n}n}{n^2}
&\equiv0\pmod p \quad\text{for $p>5$}.
\label{st2}
\end{align}
\end{lemma}

\begin{proof}
The first congruence follows from specialization $N=(p-1)/2$
of Staver's identity~\cite{St}
$$
\sum_{n=1}^N \binom{2n}{n} \frac{1}{n}
= \frac{N+1}3 \binom{2N+1}{N} \sum_{n=1}^N \frac{1}{n^2 {\binom{N}{n}}^2}.
$$
The second congruence is the modulo~$p$ reduction
of Tauraso's congruence in~\cite[Theorem~4.2]{Ta}.
It is interesting to mention that the latter follows from
the $N=p$ specialization of another combinatorial identity,
$$
\sum_{n=1}^N\binom{2n}n\frac{n^2}{4N^4+n^4}\prod_{k=1}^{n-1}\frac{N^4-k^4}{4N^4+k^4}
=\frac2{5N^2},
$$
conjectured by Borwein and Bradley~\cite{BoBr} and proved by Almkvist and Granville~\cite{AG}.
\end{proof}

Denote $q(x)=q_p(x)=(x^{p-1}-1)/p$ the Fermat quotient of $x\in\mathbb Z_p^*$.

\begin{lemma}
\label{standard2}
The following congruences is true\textup:
\begin{equation}
\sum_{n=0}^{(p-3)/2}\frac{2^{-2n}\binom{2n}n}{2n+1}
\equiv-(-1)^{(p-1)/2}q_p(2)\pmod p \quad\text{for $p>2$}.
\label{st3}
\end{equation}
\end{lemma}

\begin{proof}
Note that $(\frac12)_n\equiv(\frac12-\frac p2)_n\pmod p$ and write
the left-hand side as
\begin{equation}
\sum_{n=0}^{(p-3)/2}\frac{(\frac12)_n}{n!(2n+1)}
\equiv\sum_{n=0}^{(p-3)/2}\frac{(\frac12-\frac p2)_n}{n!(2n+1)}\pmod p.
\label{st3-1}
\end{equation}
The latter is nothing else but a terminating hypergeometric series with one term missing,
$$
\sum_{n=0}^{N-1}\frac{(-N)_n(\frac12)_n}{n!(\frac32)_n}
=-\frac{(-1)^N}{2N+1}+{}_2F_1\bigl(-N,\,\tfrac12;\,\tfrac32;\,1\bigr)
\quad\text{where $N=(p-1)/2$};
$$
it can be summed with the help of the Chu--Vandermonde theorem~\cite[Eq.~(1.7.7)]{Sl}:
\begin{equation}
\sum_{n=0}^{(p-3)/2}\frac{(\frac12-\frac p2)_n}{n!(2n+1)}
=-\frac{(-1)^{(p-1)/2}}p+\frac{(1)_{(p-1)/2}}{p(\frac12)_{(p-1)/2}}
\label{st3-2}
\end{equation}
Finally, recall Morley's congruence~\cite{Mo},
\begin{equation}
\frac{(\frac12)_{(p-1)/2}}{(1)_{(p-1)/2}}
=\binom{p-1}{\frac{p-1}2}2^{-(p-1)}
\equiv(-1)^{(p-1)/2}2^{p-1}\pmod{p^2}.
\label{st3-3}
\end{equation}
Combining \eqref{st3-1}, \eqref{st3-2}, \eqref{st3-3} and $2^{p-1}\equiv1\pmod p$
we arrive at~\eqref{st3}.
\end{proof}

\begin{lemma}
\label{pre-standard3}
For a prime $p>3$, let $x$ be a rational number such that both
$x$ and $1-x$ do not involve $p$ in their prime factorizations
and $1-x$ is a quadratic residue modulo~$p$. Take $y$ such
that $y^2\equiv1-x\pmod p$. Then
\begin{align}
\sum_{n=1}^{p-1}\binom{2n}n\biggl(\frac x4\biggr)^n
&\equiv\frac{px}{2(1-x)}\bigl(-q(x)+q(y+1)(y+1)-q(y-1)(y-1)\bigr)\pmod{p^2}.
\label{M0}
\end{align}
\end{lemma}

\begin{proof}
It is well known that
\begin{align}
\sum_{n=1}^{p-1}\frac{x^n}n
&=\sum_{n=1}^{p-1}\frac{(1)_{n-1}x^n}{(1)_n}
\equiv\sum_{n=1}^{p-1}\frac{(1-p)_{n-1}x^n}{(1)_n}\pmod p
\nonumber\\
&=-\frac1p\sum_{n=1}^{p-1}\frac{(-p)_nx^n}{(1)_n}
=\frac1p\biggl(1-x^p-\sum_{n=0}^p\binom pn(-x)^n\biggr)
\nonumber\\
&=\frac1p\bigl(1-x^p-(1-x)^p\bigr)
\label{M1}
\end{align}
and, by replacing $x$ with $-x$ and taking the appropriate linear
combination of the two expressions,
\begin{equation}
\sum_{k=1}^{(p-1)/2}\frac{x^{2k-1}}{2k-1}
\equiv\frac1{2p}\bigl((1+x)^p-(1-x)^p-2x^p\bigr)\pmod p.
\label{M2}
\end{equation}

Consider
$$
F(n,k)=\frac{(\frac12-k)_n}{(1)_n}\,\frac{x^n}{(1-x)^k}
\quad\text{and}\quad
G(n,k)=-\frac{(\frac32-k)_{n-1}}{(1)_{n-1}}\,\frac{x^n}{(1-x)^k}.
$$
They satisfy
\begin{equation}
\sum_{n=1}^{p-1}\binom{2n}n\biggl(\frac x4\biggr)^n
=\sum_{n=1}^{p-1}F(n,0)
\label{M3}
\end{equation}
and we also have
\begin{equation}
F(n,k-1)-F(n,k)=G(n+1,k)-G(n,k).
\label{00}
\end{equation}
The latter relation expresses the fact that $F(n,k)$ and $G(n,k)$ form
a \emph{WZ-pair}.

Summing \eqref{00} over $n=0,1,\dots,p-1$ we obtain
\begin{equation}
\sum_{n=0}^{p-1}F(n,k-1)-\sum_{n=0}^{p-1}F(n,k)=G(p,k)-G(0,k)=G(p,k);
\label{M4}
\end{equation}
here $G(0,k)=0$ because its expression involves $(1)_{-1}$ in the denominator.
Summing the result~\eqref{M4} over $k=1,\dots,\frac{p+1}2$ we arrive at
\begin{equation}
\sum_{n=1}^{p-1}F(n,0)=-F(0,0)+\sum_{n=0}^{p-1}F(n,\tfrac{p+1}2)
+\sum_{k=1}^{(p+1)/2}G(p,k).
\label{M5}
\end{equation}

For $n=1,2,\dots,p-1$,
$$
F(n,\tfrac{p+1}2)=\frac{(-\frac p2)_n}{(1)_n}\,\frac{x^n}{(1-x)^{(p+1)/2}}
=-\frac p2\,\frac{(1-\frac p2)_{n-1}}{(1)_n}\,\frac{x^n}{(1-x)^{(p+1)/2}}
$$
and since $(1-\frac p2)_{n-1}\equiv(1)_{n-1}\pmod p$, we obtain
\begin{equation}
\sum_{n=1}^{p-1}F(n,\tfrac{p+1}2)
\equiv-\frac p{2(1-x)^{(p+1)/2}}\sum_{n=1}^{p-1}\frac{x^n}n
\equiv-\frac{1-x^p-(1-x)^p}{2(1-x)^{(p+1)/2}}\pmod{p^2}
\label{M6}
\end{equation}
by \eqref{M1}. In addition,
\begin{equation}
F(0,0)=1 \quad\text{and}\quad
F(0,\tfrac{p+1}2)=\frac1{(1-x)^{(p+1)/2}}.
\label{M7}
\end{equation}

Furthermore, we have
$$
G(p,k)=-\frac{(\frac12-k)_p}{(\frac12-k)(p-1)!}\,\frac{x^p}{(1-x)^k}.
$$
Note that the multiples $\frac12-k+j$, $j=0,1,\dots,p-1$, runs through the complete
residue system modulo~$p$, with exactly one of them, $\frac p2$ for $j=\frac{p-1}2+k$, divisible by~$p$;
therefore $(\frac12-k)_p\equiv\frac12p!\pmod{p^2}$ for $k=1,\dots,\frac{p-1}2$.
Since $\frac12-k$ is coprime with~$p$ for those~$k$, we have
\begin{equation*}
G(p,k)\equiv\frac p{2k-1}\,\frac{x^p}{(1-x)^k}\pmod{p^2}.
\end{equation*}
If $k=\frac{p+1}2$, then
$$
\frac{(\frac12-k)_p}{\frac12-k}
=(1-\tfrac p2)_{p-1}\equiv(p-1)!\pmod{p^2},
$$
because $(1-\varepsilon)_{p-1}=(1)_{p-1}\cdot(1-\varepsilon H_{p-1}+O(\varepsilon^2))$
and the harmonic number $H_{p-1}=1+\frac12+\dots+\frac1{p-1}\equiv0\pmod p$
(see, for example, \eqref{M1} with $x=1$).
Thus, we have
\begin{equation*}
G(p,\tfrac{p+1}2)\equiv-\frac{x^p}{(1-x)^{(p+1)/2}}\pmod{p^2},
\end{equation*}
so that
\begin{align}
\sum_{k=1}^{(p+1)/2}G(p,k)
&\equiv px^p(1-x)^{-1/2}\sum_{k=1}^{(p-1)/2}\frac{\bigl((1-x)^{-1/2}\bigr)^{2k-1}}{2k-1}
-\frac{x^p}{(1-x)^{(p+1)/2}}
\nonumber\displaybreak[2]\\
&\equiv x^p\cdot\frac{\bigl(1+(1-x)^{-1/2}\bigr)^p-\bigl(1-(1-x)^{-1/2}\bigr)^p-2(1-x)^{-p/2}}{2(1-x)^{1/2}}
\nonumber\\ &\qquad
-\frac{x^p}{(1-x)^{(p+1)/2}}\pmod{p^2}
\nonumber\\
&=x^p\cdot\frac{\frac12\bigl((1-x)^{1/2}+1\bigr)^p-\frac12\bigl((1-x)^{1/2}-1\bigr)^p-2}
{(1-x)^{(p+1)/2}}
\label{M10}
\end{align}
by \eqref{M2}.

Substituting \eqref{M6}, \eqref{M7} and \eqref{M10} into \eqref{M5} we obtain,
modulo~$p^2$,
\begin{align}
\sum_{n=1}^{p-1}F(n,0)
&\equiv-1+\frac{1+(1-x)^p-3x^p+x^p\bigl((1-x)^{1/2}+1\bigr)^p-x^p\bigl((1-x)^{1/2}-1\bigr)^p}
{2(1-x)^{(p+1)/2}}
\nonumber\\
&=\frac{\bigl((1-x)^{(p-1)/2}-1\bigr)^2}{2(1-x)^{(p-1)/2}}
\nonumber\\ &\;
+\frac x{2(1-x)}\,
\frac{1-3x^{p-1}+x^{p-1}\bigl((1-x)^{1/2}+1\bigr)^p-x^{p-1}\bigl((1-x)^{1/2}-1\bigr)^p}{(1-x)^{(p-1)/2}}
\nonumber\\
&=\frac{(y^{p-1}-1)^2}{2y^{p-1}}
+\frac x{2(1-x)}\,\frac{1-3x^{p-1}+x^{p-1}(y+1)^p-x^{p-1}(y-1)^p}{y^{p-1}}.
\label{M11}
\end{align}
Noting that $(y^{p-1}-1)^2=p^2q(y)\equiv0\pmod{p^2}$,
\begin{align*}
&
\frac{1-3x^{p-1}+x^{p-1}(y+1)^p-x^{p-1}(y-1)^p}p
\\ &\quad
=-\frac{x^{p-1}-1}p(3-(y+1)^p+(y-1)^p)
\\ &\quad\qquad
+\frac{(y+1)^{p-1}-1}p(y+1)
-\frac{(y-1)^{p-1}-1}p(y-1)
\\ &\quad
\equiv-q(x)(3-(y+1)+(y-1))+q(y+1)(y+1)-q(y-1)(y-1)\pmod p
\\ &\quad
=-q(x)+q(y+1)(y+1)-q(y-1)(y-1),
\end{align*}
and $y^{p-1}\equiv1\pmod p$, we obtain the required congruence~\eqref{M0}
from \eqref{M3} and~\eqref{M11}.
\end{proof}

\begin{lemma}
\label{standard3}
The following congruences are valid\textup:
\begin{align}
\sum_{n=1}^{p-1}\frac{(-1)^n2^n\binom{2n}n}n
&\equiv-4q_p(2)\pmod p \quad\text{for $p>2$},
\label{st4}
\\
3\sum_{n=1}^{p-1}(-1)^n2^n\binom{2n}n
&\equiv-4q_p(2)\pmod{p^2} \quad\text{for $p>2$}.
\label{st5}
\end{align}
\end{lemma}

\begin{proof}
It is shown in~\cite[Theorem~1.2]{ST} that for $m$ in~$\mathbb Z_p^*$,
\begin{equation}
\sum_{n=1}^{p-1}\frac{(-1)^n\binom{2n}n}{nm^n}
\equiv\frac2m\cdot\frac{m^p-V_p(m)}p\pmod p,
\label{st4-1}
\end{equation}
where the sequence $V_k(x)$ is defined by $V_0(x)=2$, $V_1(x)=x$,
and $V_k(x)=x(V_{k-1}(x)+V_{k-2}(x))$ for $k\ge2$. Although the
theorem is stated for $m\in\mathbb Z$ only, the proof does not make
use of this integrality: we can apply it for $m=1/2$ as well.
In this case $V_k(1/2)=1+(-1)^k/2^k$, so that the right-hand side
of~\eqref{st4-1} becomes~\eqref{st4} if we additionally use $2^{p-1}\equiv1\pmod p$.

The congruence~\eqref{st5} is clear for $p=3$, while for $p>3$ it follows
from specialization $x=-8$, $y=3$ of \eqref{M0} and noting that
\begin{align*}
q_p(-8)
&=\frac{2^{3(p-1)}-1}p
=\frac{2^{p-1}-1}p\,(2^{2(p-1)}+2^{p-1}+1)
\equiv 3q_p(2)\pmod p,
\\
q_p(4)
&=\frac{2^{2(p-1)}-1}p
=\frac{2^{p-1}-1}p\,(2^{p-1}+1)
\equiv 2q_p(2)\pmod p.
\end{align*}
\vskip-\baselineskip
\end{proof}

\section{Proofs of the supercongruences}
\label{s3}

\begin{proof}[Proof of~\textup{\eqref{J1}}]
Take
$$
F(n,k)=(3n+2k+1)\frac{(\frac12)_n(\frac12+k)_n^2}{(1)_n^3}2^{2n}
\quad\text{and}\quad
G(n,k)=-\frac{(\frac12)_n(\frac12+k)_{n-1}^2}{(1)_{n-1}^3}2^{2n}.
$$
Then we have
\begin{equation}
\sum_{n=0}^{(p-1)/2}\frac{(\frac12)_n^3}{n!^3}(3n+1)2^{2n}
=\sum_{n=0}^{(p-1)/2}F(n,0)
\label{S1}
\end{equation}
and \eqref{00}, so that $F(n,k)$ and $G(n,k)$ form a WZ-pair.
Summing \eqref{00} over $n=0,1,\dots,\frac{p-1}2$, we obtain
\begin{equation}
\sum_{n=0}^{(p-1)/2}F(n,k-1)
-\sum_{n=0}^{(p-1)/2}F(n,k)
=G(\tfrac{p+1}2,k)-G(0,k)
=G(\tfrac{p+1}2,k).
\label{01}
\end{equation}
Furthermore, for $k=1,2,\dots,\frac{p-1}2$ we have
\begin{equation*}
G(\tfrac{p+1}2,k)
=-\frac{(\frac12)_{(p-1)/2}(\frac12+k)_{(p-1)/2}^2}{(1)_{(p-1)/2}^3}2^{p+1}
\equiv 0\pmod{p^3},
\end{equation*}
because each of the three Pochhammer products in the numerator is divisible by~$p$
while the denominator, $\bigl(\frac{p-1}2\bigr)!^3$, is coprime with~$p$.
Comparing this result with~\eqref{01}, as in the proof of Theorem~1 in~\cite{Zu}, we see that
\begin{equation*}
\sum_{n=0}^{(p-1)/2}F(n,0)
\equiv\sum_{n=0}^{(p-1)/2}F(n,1)
\equiv\sum_{n=0}^{(p-1)/2}F(n,2)
\equiv\dots
\equiv\sum_{n=0}^{(p-1)/2}F(n,\tfrac{p-1}2)\pmod{p^3},
\end{equation*}
hence we can replace, modulo $p^3$, our sum~\eqref{S1} by
\begin{align}
\sum_{n=0}^{(p-1)/2}F(n,\tfrac{p-1}2)
&=\sum_{n=0}^{(p-1)/2}(3n+p)\frac{(\frac12)_n(\frac p2)_n^2}{(1)_n^3}2^{2n}
\nonumber\displaybreak[2]\\
&=p+p\sum_{n=1}^{(p-1)/2}\frac{(\frac12)_n(\frac p2)_n^2}{(1)_n^3}2^{2n}
+3\sum_{n=1}^{(p-1)/2}n\frac{(\frac12)_n(\frac p2)_n^2}{(1)_n^3}2^{2n}
\nonumber\displaybreak[2]\\
&=p+\frac{p^3}4\sum_{n=1}^{(p-1)/2}\frac{(\frac12)_n(1+\frac p2)_{n-1}^2}{(1)_n^3}2^{2n}
+\frac{3p^2}4\sum_{n=1}^{(p-1)/2}n\frac{(\frac12)_n(1+\frac p2)_{n-1}^2}{(1)_n^3}2^{2n}.
\label{03}
\end{align}
(Note that, in contrast with the proofs in~\cite{Zu}, the newer sum is not reduced
to a single term.) Comparing the resulted expression~\eqref{03} for~\eqref{S1}
we see that \eqref{J1} is equivalent to
\begin{equation}
\sum_{n=1}^{(p-1)/2}n\frac{(\frac12)_n(1+\frac p2)_{n-1}^2}{(1)_n^3}2^{2n}
\equiv0\pmod p \quad\text{for $p>3$}.
\label{04}
\end{equation}
On noting that $(1+\frac p2)_{n-1}\equiv(1)_{n-1}\pmod p$, we reduce \eqref{04}
to its equivalent
\begin{equation*}
\sum_{n=1}^{(p-1)/2}\frac{(\frac12)_n}{n(1)_n}2^{2n}
=\sum_{n=1}^{(p-1)/2}\frac{\binom{2n}n}n
\equiv0\pmod p \quad\text{for $p>3$},
\end{equation*}
which is exactly~\eqref{st1}.
\end{proof}

\begin{proof}[Proof of~\textup{\eqref{J2}}]
The proof is very similar. This time we take
\begin{gather*}
F(n,k)=(10n^2+12nk+4k^2+6n+4k+1)\frac{(\frac12)_n(\frac12+k)_n^4}{(1)_n^5}(-1)^n2^{2n}
\quad\text{and}\\
G(n,k)=(n+2k-1)\frac{(\frac12)_n(\frac12+k)_{n-1}^4}{(1)_{n-1}^5}(-1)^n2^{2n+1}
\end{gather*}
with the motive
\begin{equation*}
\sum_{n=0}^{(p-1)/2}\frac{(\frac12)_n^5}{n!^5}(10n^2+6n+1)(-1)^n2^{2n}
=\sum_{n=0}^{(p-1)/2}F(n,0)
\end{equation*}
and~\eqref{00}. Then, as above, we find that
\begin{align}
\sum_{n=0}^{(p-1)/2}F(n,0)
&\equiv\sum_{n=0}^{(p-1)/2}F(n,\tfrac{p-1}2)\pmod{p^5}
\nonumber\displaybreak[2]\\
&=\sum_{n=0}^{(p-1)/2}(10n^2+6np+p^2)\frac{(\frac12)_n(\frac p2)_n^4}{(1)_n^5}(-1)^n2^{2n}
\nonumber\displaybreak[2]\\
&=p^2+\frac{p^6}{16}\sum_{n=1}^{(p-1)/2}\frac{(\frac12)_n(1+\frac p2)_{n-1}^4}{(1)_n^5}(-1)^n2^{2n}
\nonumber\\ &\qquad
+\frac{3p^5}8\sum_{n=1}^{(p-1)/2}n\frac{(\frac12)_n(1+\frac p2)_{n-1}^4}{(1)_n^5}(-1)^n2^{2n}
\nonumber\\ &\qquad
+\frac{5p^4}8\sum_{n=1}^{(p-1)/2}n^2\frac{(\frac12)_n(1+\frac p2)_{n-1}^4}{(1)_n^5}(-1)^n2^{2n}
\nonumber
\end{align}
and our task is to show that
\begin{equation*}
\sum_{n=1}^{(p-1)/2}n^2\frac{(\frac12)_n(1+\frac p2)_{n-1}^4}{(1)_n^5}(-1)^n2^{2n}
\equiv0\pmod p \quad\text{for $p>5$}.
\end{equation*}
Using $(1+\frac p2)_{n-1}\equiv(1)_{n-1}\pmod p$ the latter reduces to
\begin{equation*}
\sum_{n=1}^{(p-1)/2}\frac{(\frac12)_n}{n^2(1)_n}(-1)^n2^{2n}
=\sum_{n=1}^{(p-1)/2}\frac{(-1)^n\binom{2n}n}{n^2}
\equiv0\pmod p \quad\text{for $p>5$},
\end{equation*}
which is \eqref{st2}.
\end{proof}

\begin{proof}[Proof of~\textup{\eqref{J4}}]
Take the WZ-pair
\begin{gather*}
F(n,k)=(3n+2k+1)\frac{(\frac12)_n(\frac12+k)_n^2(\frac12)_k}{(1)_n^2(1+2k)_n(1)_k}(-1)^n2^{3n}
\quad\text{and}\\
G(n,k)=\frac{(\frac12)_n(\frac12+k)_{n-1}^2(\frac12)_k}{(1)_{n-1}^2(1+2k)_{n-1}(1)_k}(-1)^n2^{3n-2}.
\end{gather*}
Summing \eqref{00} over $n=0,1,\dots,p-1$ we get
\begin{equation}
\sum_{n=0}^{p-1}F(n,k-1)
-\sum_{n=0}^{p-1}F(n,k)
=G(p,k)-G(0,k)
=G(p,k).
\label{21}
\end{equation}
Summing \eqref{21} further over $k=1,2,\dots,\frac{p-1}2$ we obtain
\begin{equation}
\sum_{n=0}^{p-1}F(n,0)
=\sum_{n=0}^{p-1}F(n,\tfrac{p-1}2)
+\sum_{k=1}^{(p-1)/2}G(p,k).
\label{22}
\end{equation}

Note that
\begin{align*}
G(p,k)
&=\frac{(\frac12+p)_{k-1}^2(1)_{2k}}{(\frac12)_k(1+p)_{2k-1}(1)_k}
\cdot\frac{(\frac12)_p^3}{p(1)_{p-1}^3}(-1)^p2^{3p-2}
\\
&=-\frac{(\frac12+p)_{k-1}^2\cdot2^{2k}}{(1+p)_{2k-1}}
\cdot\frac{(\frac12)_p^3}{p(p-1)!^3}2^{3p-2}.
\end{align*}
For $k=1,2,\dots,\frac{p-1}2$ neither $(\frac12+p)_{k-1}$ nor $(1+p)_{2k-1}$
is divisible by $p$; in addition, we have
$(\frac12+p)_{k-1}\equiv(\frac12)_{k-1}\pmod p$ and $(1+p)_{2k-1}\equiv(1)_{2k-1}\pmod p$.
Since the factor in $G(p,k)$ independent of~$k$ is divisible by~$p^2$, we conclude that
\begin{align*}
\sum_{k=1}^{(p-1)/2}G(p,k)
&\equiv-\frac{(\frac12)_p^3}{p(p-1)!^3}2^{3p-2}
\sum_{k=1}^{(p-1)/2}\frac{(\frac12)_{k-1}^2}{(2k-1)!}2^{2k}
\\
&\equiv-p^2\binom{2p-1}{p-1}2^{-3p+3}
\sum_{k=1}^{(p-1)/2}\frac{2^{-2(k-1)}\binom{2k-2}{k-1}}{2k-1}\pmod{p^3}.
\end{align*}
Finally, we use $\binom{2p-1}{p-1}\equiv1\pmod p$, $2^{p-1}\equiv1\pmod p$ and
the congruence~\eqref{st3} to get
\begin{equation}
\sum_{k=1}^{(p-1)/2}G(p,k)
\equiv(-1)^{(p-1)/2}p(2^{p-1}-1)\pmod{p^3}.
\label{23}
\end{equation}

Furthermore,
\begin{equation}
F(0,\tfrac{p-1}2)
=p\frac{(\frac12)_{(p-1)/2}}{(1)_{(p-1)/2}}
\equiv (-1)^{(p-1)/2}p2^{p-1}\pmod{p^3}
\label{24}
\end{equation}
by~\eqref{st3-3} and
\begin{align*}
F(n,\tfrac{p-1}2)
&=(3n+p)\frac{(\frac12)_n(\frac p2)_n^2(\frac12)_{(p-1)/2}}{(1)_n^2(p)_n(1)_{(p-1)/2}}(-1)^n2^{3n}
\\
&\equiv(3n+p)\frac{(\frac12)_n(1+\frac p2)_{n-1}^2}{(1)_n^2(1+p)_{n-1}}(-1)^n2^{3n-2}
\cdot(-1)^{(p-1)/2}p2^{p-1}\pmod{p^3}.
\end{align*}
Using $(1+\varepsilon)_{n-1}=(1)_{n-1}\cdot(1+\varepsilon H_{n-1}+O(\varepsilon^2))$,
where $H_{n-1}=1+\frac12+\dots+\frac1{n-1}$, with $\varepsilon=\frac p2$ and~$p$ we find that
$$
\frac{(1+\frac p2)_{n-1}^2}{(1+p)_{n-1}}
=\frac{(1)_{n-1}^2\cdot(1+\frac p2H_{n-1}+O(p^2))^2}{(1)_{n-1}\cdot(1+pH_{n-1}+O(p^2))}
=(1)_{n-1}\cdot(1+O(p^2))
\equiv(1)_{n-1}\pmod{p^2},
$$
hence
\begin{align*}
F(n,\tfrac{p-1}2)
&\equiv(3n+p)\frac{(\frac12)_n}{(1)_nn}(-1)^n2^{3n-2}
\cdot(-1)^{(p-1)/2}p2^{p-1}\pmod{p^3}
\\
&=\frac14(3n+p)\frac{(-1)^n2^n\binom{2n}n}n
\cdot(-1)^{(p-1)/2}p2^{p-1}
\end{align*}
and
\begin{align}
\sum_{n=1}^{p-1}F(n,\tfrac{p-1}2)
&\equiv(-1)^{(p-1)/2}p2^{p-1}\biggl(\frac34\sum_{n=1}^{p-1}(-1)^n2^n\binom{2n}n
+\frac p4\sum_{n=1}^{p-1}\frac{(-1)^n2^n\binom{2n}n}n\biggr)
\nonumber\\
&\equiv(-1)^{(p-1)/2}p2^{p-1}\cdot2(1-2^{p-1})
\equiv(-1)^{(p-1)/2}p\cdot2(1-2^{p-1})\pmod{p^3}
\label{25}
\end{align}
by \eqref{st4}, \eqref{st5} and $2^{p-1}\equiv1\pmod p$.

Substituting \eqref{23}, \eqref{24} and  \eqref{25} into~\eqref{22}
we obtain
\begin{align*}
\sum_{n=0}^{p-1}F(n,0)
&\equiv(-1)^{(p-1)/2}p\bigl(2^{p-1}+2(1-2^{p-1})+(2^{p-1}-1)\bigr)\pmod{p^3}
\\
&=(-1)^{(p-1)/2}p,
\end{align*}
the required congruence.
\end{proof}

\section{``Divergent'' Ramanujan-type series}
\label{sdiv}

In~\cite{Zu}, the second-named author generalized an observation of L.~Van Hamme
about Ramanujan-type identities for $1/\pi$ and $1/\pi^2$. The idea is to associate
with each such identity
\begin{equation*}
\sum_{n=0}^\infty A_n(a+bn)z^n=\frac{r\sqrt{d}}{\pi}
\quad\text{or}\quad
\sum_{n=0}^\infty A_n(a+bn+cn^2)z^n=\frac{r\sqrt{d}}{\pi^2},
\end{equation*}
where $a$, $b$, $c$, $z$, and $r$ are rational and $A_n$ is a related Pochhammer
ratio (or, more generally, an Ap\'ery-like sequence; cf.~\cite{Zu}), the supercongruence
\begin{equation*}
\sum_{n=0}^{p-1}A_n(a+bn)z^n\overset?\equiv a\biggl(\frac{-d}{p}\biggr)p \pmod{p^3}
\end{equation*}
or
\begin{equation}
\sum_{n=0}^{p-1}A_n(a+bn+cn^2)z^n\overset?\equiv a\biggl(\frac{d}{p}\biggr)p^2 \pmod{p^5},
\end{equation}
respectively, for all $p\ge p_0$. Recently \cite{Gu2}, the first-named author went even further
and considerably extended the pattern; however this remains an unproven observation.

The general machinery for proving Ramanujan-like series for $1/\pi$
\cite{BoBo,CCL,Zu1} produces, in several cases, \emph{divergent instances} like
\begin{equation}
\sum_{n=0}^{\infty}\frac{(\frac12)_n^3}{(1)_n^3}(3n+1)2^{2n}\,\text{``$=$''}\,\frac{-2i}{\pi},
\qquad
\sum_{n=0}^{\infty}\frac{(\frac12)_n^3}{(1)_n^3}(3n+1)(-1)^n2^{3n}\,\text{``$=$''}\,\frac1\pi.
\label{3F2-div3}
\end{equation}
The summations in~\eqref{3F2-div3} have to be understood as the analytic continuation
of the corresponding $_3F_2$-hypergeometric series; for example, the second formula in~\eqref{3F2-div3}
can be given by
$$
\frac1{2\pi i}\int_{-1/4-i\infty}^{-1/4+i\infty}
\frac{(\frac12)_s^3}{(1)_s^2}\,\Gamma(-s)(3s+1)2^{3s}\,\mathrm ds
=\frac1\pi.
$$
The first appearance of divergent series for $1/\pi$ is \cite[p.~371]{BoBo}.
In view of the observation from~\cite{Zu},
the formulae in~\eqref{3F2-div3} motivate our supercongruences~\eqref{J1a}
and~\eqref{3F2-zu3}, respectively.

Curiously, our study of the first identity in~\eqref{3F2-div3} led us to
``complex'' convergent Ramanujan-type formulae for $1/\pi$ such as
\begin{equation}
\sum_{n=0}^{\infty}\frac{(\frac12)_n^3}{(1)_n^3}
\biggl(\frac{105-21\sqrt{-7}}{32}n+\frac{49-13\sqrt{-7}}{64}\biggr)
\cdot\biggl(\frac{47+45\sqrt{-7}}{128}\biggr)^n
=\frac{\sqrt{7}}{\pi}.
\label{3F2-divm}
\end{equation}
As far as we know such formulae do not exist in the literature.
Note that application of the quadratic transformation
\begin{equation*}
{}_3F_2\biggl(\begin{matrix}
\frac12, \, \frac12, \, \frac12 \\
1, \, 1 \end{matrix} \biggm| z\biggr)
=(1-z)^{-1/2}\cdot{}_3F_2\biggl(\begin{matrix}
\frac14, \, \frac12, \, \frac34 \\
1, \, 1 \end{matrix} \biggm|\frac{-4z}{(1-z)^2}\biggr)
\end{equation*}
(the method used in~\cite{CZ} and~\cite{Zu0}) translates~\eqref{3F2-divm} into
the identity
\begin{equation}
\sum_{n=0}^\infty \frac{(\frac12)_n(\frac14)_n(\frac34)_n}{(1)_n^3}
(35n+8)\biggl(\frac43\biggr)^{4n}\,\text{``$=$''}\,-\frac{18i}{\pi}
\label{3F2-divn}
\end{equation}
which has to be understood as the analytic continuation of the hypergeometric
series on the left-hand side to $\mathbb C\setminus(-\infty,0]$
and which serves as the prototype of~\eqref{3F2-zun}.
All identities for $1/\pi$, including the divergent and complex instances
\eqref{3F2-div3}--\eqref{3F2-divn} and others, can be proven by the
modular argument. We plan to address these issues in another project.
In~\cite{Gu1}, the first-named author gives proofs of several divergent hypergeometric
formulae for $1/\pi$ and $1/\pi^2$ using a version of the Wilf--Zeilberger algorithm.

The theory developed in \cite{Gu}, allows us to obtain numerically
the parameters of the divergent Ramanujan-like series for $1/\pi^2$ as well.
For example, the expansion
\begin{equation*}
\sum_{n=0}^{\infty}\frac{(\frac12)_{n+x}^5}{(1)_{n+x}^5}(a+b(n+x)+c(n+x)^2)(-1)^nz^{n+x}
=\frac1{\pi^2}-\frac k2x^2+O(x^4)
\quad\text{as $x\to0$},
\end{equation*}
where $z>0$, corresponds to the case $s=t=1/2$ and $u=-1$ of~\cite[Exp.~1.2]{Gu}, and
defines $z$, $a$, $b$ and $c$ as functions of~$k$. All these quantities admit
a natural parametrization by means of~$\tau$ where $\tau^2=c^2/(1+z)$;
see \cite[Eq.~(3.47)]{Gu} for details.
For this case, the first-named author discovered experimentally that the relation
$z(\tau_1)z(\tau_2)=1$ implies
\begin{equation}
(k_1+1)\tau_2=(k_2+1)\tau_1, \quad (k_1+1)(k_2+1)+8=4 \tau_1 \tau_2,
\label{ktau}
\end{equation}
and also
\begin{equation}
c(\tau_2)=\frac{\tau_2}{\tau_1} \cdot \frac{c}{\sqrt{z}}(\tau_1), \quad
b(\tau_2)=\frac{\tau_2}{\tau_1} \cdot \frac{c-b}{\sqrt{z}}(\tau_1), \quad
a(\tau_2)=\frac{\tau_2}{\tau_1} \cdot \frac{c-2b+4a}{4 \sqrt{z}}(\tau_1).
\label{par-abc}
\end{equation}
The choice $\tau_1=\sqrt{5}$ corresponds to the Ramanujan-like series
\begin{equation*}
\sum_{n=0}^{\infty}\frac{(\frac12)_n^5}{(1)_n^5}(20n^2+8n+1)\frac{(-1)^n}{2^{2n}}=\frac{8}{\pi^2}
\end{equation*}
proven in~\cite{Gu0}, in which case $k=1$ and $z(\tau_1)=1/4$.
This series suggests the existence of a ``divergent'' companion
with $z(\tau_2)=4$, $\tau_2=\sqrt{5}/2$, $k=0$, $c=5/2$, $b=3/2$ and $a=1/4$;
the values $\tau$ and $k$ are found from \eqref{ktau} and $c$, $b$ and $a$ from~\eqref{par-abc}.
The resulting set corresponds to the series
\begin{equation*}
\sum_{n=0}^{\infty} \frac{(\frac12)_n^5}{(1)_n^5}(10n^2+6n+1)(-1)^n2^{2n}\,\text{``$=$''}\,\frac{4}{\pi^2}
\end{equation*}
with the left-hand side understood as the analytic continuation of the participating
hypergeometric series to $\mathbb C\setminus[1,+\infty)$. A similar duality for the
$_3F_2$-evaluations of $1/\pi$ can be explained by the modular origin of the corresponding
hypergeometric series, like the one we give for~\eqref{3F2-div3}. The duality mechanism for
the $_5F_4$-examples remains a mystery.

As already pointed out in~\cite{Zu}, \emph{all} Ramanujan-type series for $1/\pi$ and their generalizations
possess unexpectedly strong arithmetic properties. In particular, these are reflected by
the supercongruences for truncated sums\,---\,it is probably not surprising to see the
examples \eqref{J1a}--\eqref{5F4-zu4}. What is more remarkable, the $p$-analogues
make no difference of their origin: whether they come from convergent or divergent
formulae. This kind of democracy as well as an apparent simplicity of the supercongruences
make them an attractive object for further investigation.

\end{document}